\documentclass[12pt]{article} 
\usepackage{amsmath}
\usepackage{amssymb}
\usepackage{amsthm}
\usepackage{mathtools}
\DeclarePairedDelimiter\ceil{\lceil}{\rceil}

\usepackage[utf8]{inputenc}
\usepackage[english]{babel}
\usepackage{color}

\newcommand{\Mod}[1]{\ (\mathrm{mod}\ #1)}

\usepackage{wrapfig}

\newtheorem{theorem}{Theorem}[section]
\newtheorem{corollary}[theorem]{Corollary}

\newtheorem{lemma}[theorem]{Lemma}
\newtheorem{definition}[theorem]{Definition}
\newtheorem{proposition}[theorem]{Proposition}

\newtheorem{claim}[theorem]{Claim}

\title{Bounds on point configurations determined by distances and dot products}
\date{\today}
\author{Slade Gunter, Eyvi Palsson, Ben Rhodes, and Steven Senger}

\begin{document}
\maketitle
\begin{abstract}
We study a family of variants of Erd\H os' unit distance problem, concerning distances and dot products between pairs of points chosen from a large finite point set. Specifically, given a large finite set of $n$ points $E$, we look for bounds on how many subsets of $k$ points satisfy a set of relationships between point pairs based on distances or dot products. We survey some of the recent work in the area and present several new, more general families of bounds.
\end{abstract}

\section{Introduction}
Paul Erd\H os introduced the single distance problem in \cite{Erd46}. It asks how often a single distance can occur between pairs of points chosen from a large finite point set in the plane. While there have been many generalizations and variants of this original problem involving distances, there has also been work on pursuing similar questions for dot products, as in \cite{GIS, Steinerberger}. Moreover, configurations involving more than two points have also seen attention. See \cite{Bahls, Fickus, IS} for examples and applications. The book by Brass, Moser, and Pach, \cite{BMP}, details many related problems on point configurations.

We offer some upper and lower bounds on general families of point configurations determined by distances or dot products between pairs of points. In what follows, if two quantities, $X(n)$ and $Y(n)$, vary with respect to some natural number parameter, $n$, then we write $X(n) \lesssim Y(n)$ if there exist constants, $C$ and $N$, both independent of $n$, such that for all $n> N$, we have $X(n)\leq CY(n)$. If $X(n) \lesssim Y(n)$ and $Y(n) \lesssim X(n)$, we write $X(n) \approx Y(n).$ By convention, we will always assume that the parameters associated to the configurations, $h$ and $k$, are like constants compared to the size parameter $n.$ 

\subsection{Background}
We begin by recording some of the best estimates for some often studied point configurations. In 1984, Spencer, Szemer\'edi, and Trotter gave the least known upper bound on the number of times that a distance can occur in a large finite point set, in \cite{SST84}.
\begin{theorem}\label{SST}[\cite{SST84}]
Given a large finite set of $n$ points in the plane, the maximum number of times that any distance can occur is bounded above by $\lesssim n^\frac{4}{3}.$
\end{theorem}
There is still a wide gap between this bound and Erd\H os' conjectured upper bound of $n^{1+\epsilon},$ for any $\epsilon>0$. However, an easy consequence of the celebrated Szemer\'edi-Trotter point-line incidence theorem from \cite{ST83} yields the same bound for dot products, which happens to be sharp. Though they did not explicitly state this bound, it is a direct corollary of their main result. We make this connection explicit below.
\begin{theorem}\label{STdp}[\cite{ST83}]
Given a large finite set of $n$ points in the plane, the maximum number of times that any nonzero dot product can occur is bounded above by $\lesssim n^\frac{4}{3}.$ Moreover, this bound is sharp.
\end{theorem}

\subsubsection{$k$-chains}
The first type of generalization that we introduce here is a {\it $k$-chain}, which is a sequence of distinct points restricted by the values of the distances or dot products determined by successive pairs. We borrow notation from related problems on distances in \cite{BIT, OT, PSS}. Specifically, if we fix a $k$-tuple of real numbers, $(\alpha_1, \alpha_2, \dots, \alpha_k)$, then a distance {\it $k$-chain} of that type is a $(k+1)$-tuple of points, $(x_1, x_2, \dots, x_{k+1}),$ such that for all $j=1,\dots, k,$ we have $|x_j- x_{j+1}|=\alpha_j.$ A dot product $k$-chain is defined similarly, except with $\alpha_j=x_j\cdot x_{j+1}$ in place of the distance. For example, if we fix a triple of real numbers, $(\alpha, \beta, \gamma)$, a distance 3-chain of that type will be a set of four points, where the distance between the first two points is $\alpha$, the distance between the middle two points is $\beta$, and the distance between last two points is $\gamma.$ We follow convention and refer to 2-chains as {\it hinges}. In \cite{BS} and \cite{S11}, Dan Barker and the fourth listed author gave the following tight bounds on the number of hinges (2-chains) in a large finite point set in the plane.

\begin{theorem}\label{hinges}[\cite{BS,S11}]
Given a large, finite set of $n$ points in $\mathbb R^2$, and a pair of nonzero real numbers $(\alpha_1, \alpha_2)$, the maximum number distance or dot product hinges of type $(\alpha_1, \alpha_2)$ is no more than $\lesssim n^2.$ Moreover, this bound is tight.
\end{theorem}

While distances and dot products have had similar behavior in these first few estimates, their paths diverge for longer chains. In \cite{PSS}, Adam Sheffer, and the second and fourth listed authors gave upper and lower bounds on the maximum possible number of occurrences of $k$-chains of a given type in large finite point sets in the plane. Many of these bounds were improved by Frankl and Kupavskii, in \cite{FK}, where they give the best known upper bounds on longer distance $k$-chains, which we record below.

\begin{theorem}\label{FKthm}[Theorem 2 in \cite{FK}]
Given a large, finite set of $n$ points in $\mathbb R^2,$ any $\epsilon >0,$ and a natural number $k$, the maximum number of distance $k$-chains of the type $(\alpha_1,\dots, \alpha_k)$ that can exist in the set is no more than
$$\lesssim
\begin{cases}
n^\frac{k+3}{3} \qquad & \text{ if } k \equiv 0 \Mod{3}, \\
n^{\frac{k+3}{3}+\epsilon} & \text{ if } k \equiv 1 \Mod{3}, \\
n^\frac{k+4}{3} & \text{ if } k \equiv 2 \Mod{3}.
\end{cases}$$
\end{theorem}

In the case of dot products, the best known upper bounds are from Shelby Kilmer, Caleb Marshall, and the fourth listed author, in \cite{KMS}. Here and in what follows, we will bound results on dot products that are nonzero, as there are degeneracies that can occur for these types of estimates when we consider zero dot products. See \cite{KMS} for more details.

\begin{theorem}\label{2chainzgeneral}[\cite{KMS}]
Given a large, finite set of $n$ points in $\mathbb R^2$ and a natural number $k$, the maximum number of dot product $k$-chains of the type $(\alpha_1,\dots, \alpha_k)$ that can exist in the set is $\lesssim n^\frac{2(k+1)}{3}.$
\end{theorem}

Although the upper bounds are clearly larger in the case of dot products, it is unclear if they are closer to or further from the truth than what is known about distances, as they exhibit fundamentally distinct behaviors for $k\geq 3,$ as the following results will show. We now turn our attention to lower bounds on the maximum number of $k$-chains of a given type we can construct.

\begin{theorem}\label{sharpFKthm}[Theorem 2 in \cite{FK}]
Given a large, finite natural number $n$, and a small natural number $k$, there exists a set of $n$ points in $\mathbb R^2,$ and a $k$-tuple $(\alpha_1,\dots,\alpha_k),$ such that the set has a number $k$-chains of the type $(\alpha_1,\dots, \alpha_k)$ that is at least
$$\gtrsim
\begin{cases}
n^\frac{k+3}{3} \qquad & \text{ if } k \equiv 0 \Mod{3}, \\
n^{\frac{k+2}{3}+\epsilon} & \text{ if } k \equiv 1 \Mod{3}, \\
n^\frac{k+4}{3} & \text{ if } k \equiv 2 \Mod{3}.
\end{cases}$$
\end{theorem}

In the case of dot products, there is a family of constructions demonstrating a lower bound for the maximum number of dot product $k$-chains in a large finite set of points in the plane.

\begin{proposition}\label{sharpDpChains}[\cite{KMS}]
Given a large, finite natural number $n$, and a small natural number $k$, there exists a set of $n$ points in $\mathbb{R}^2$ and a $k$-tuple of nonzero real numbers $(\alpha_1,\dots, \alpha_k)$, for which there are at least $n^{\ceil{(k+1)/2}}$ instances of $k$-chains of the type $(\alpha_1,\dots, \alpha_k)$ in the set.
\end{proposition}

Comparing these results shows that the behaviors of distances and dot products are distinct for $k\geq 4.$ Indeed, if we fix $k=4,$ Theorem \ref{sharpDpChains} gives that there exists a set of $n$ points in the plane with $\gtrsim n^3$ occurrences of some dot product 4-chain, while Theorem \ref{FKthm} shows that the number of occurrences of any type of distance 4-chain can never exceed $\lesssim n^\frac{7}{3}+\epsilon,$ for any $\epsilon>0.$ See also \cite{CS}, \cite{IR}, and \cite{HIKR} for some related results in vector spaces over finite fields and modules over finite rings of integers.

So far, we have only considered sequences of points. However, to deal with more general relationships between points we will need a way to organize our information.

\subsubsection{Graphs}
Let a {\it graph} $G$ be defined by a pair of sets, $V(G)$ and $E(G)$, where the first set has elements called {\it vertices}, and the second set consists of two-element subsets of $V(G)$ called {\it edges}. The {\it degree} of a vertex is the number of edges it is an element of. In a given graph, if two vertices are listed in the same edge, we will call them {\it adjacent}. Finally, we will define an {\it edge-weighted graph} as a graph $G$ paired with a {\it weight function}, $w:E(G)\rightarrow \mathbb R.$ Using edge-weighted graphs, we can specify a $k$-chain as a graph on $k+1$ vertices with edges connecting successive vertices. Each edge will have a weight corresponding to the $\alpha_j$. A {\it subgraph} of $G$ is a graph $G'$ with $V(G')\subseteq V(G)$ and $E(G')\subseteq E(G).$ Further, $G'$ is an {\it induced subgraph} if it is a subgraph with all possible edges from $E(G)$ for the vertices present in $V(G').$

\subsubsection{Three dimensions}

Thus far, we have only considered point sets in the plane. The reason for this is that in dimensions four and higher, there is a simple construction due to Lenz that renders the unit distance problem trivial as stated. There is a similar such construction in three dimensions for dot products. These constructions can then be modified to quickly make many related questions trivial unless there are further restrictions on the point sets. See \cite{BMP} or \cite{KMS} for more details. That said, we include here a few distance results in three dimensions that do not appear to follow from such simple adaptations of the Lenz example.

Let $u_3(n)$ denote the maximum number of times the unit distance can occur in any large finite set of $n$ points in $\mathbb R^3.$ In \cite{Erd60}, Erd\H os proved that for any large $n$, $u_3(n)\gtrsim n^\frac{4}{3} \log \log n$ pairs of points separated by the unit distance. In \cite{Zahl}, Josh Zahl showed $u_3(n)\lesssim n^{\frac{295}{197}+\epsilon},$ for any $\epsilon >0.$

Currently, the best known upper bounds on the number of occurrences of a given type of $k$-chain in a large finite set of $n$ points in $\mathbb R^3$ are due to Frankl and Kupavskii, in \cite{FK}. 
\begin{theorem}\label{FKthm3}[\cite{FK}]
Given a large, finite set of $n$ points in $\mathbb R^3$ and a natural number $k$, the maximum number of distance $k$-chains of the type $(\alpha_1,\dots, \alpha_k)$ that can exist in the set, is $\lesssim n^{\frac{k}{2}+1+\epsilon},$ for any $\epsilon >0.$
\end{theorem}
For even $k$, this essentially settles the problem, as it nearly matches the lower bound $\gtrsim n^{\frac{k}{2}+1},$ from \cite{PSS}. In the case that $k$ is odd, there is a wider gap between the upper and lower bounds, but Frankl and Kupavskii have the best results. We refer the interested reader to Theorem 4 and Proposition 2 in \cite{FK} for full details, as the lower bound is somewhat difficult to express concisely.




\subsection{Main results}

We now introduce the definitions we use to quantify many of the results in this paper. Given an edge-weighted graph $G$, we say that a $k$-tuple of points $X=(x_1, \dots, x_k)$ forms a {\it distance $G$-configuration} if there is a bijection from $\varphi:V(G)\rightarrow X$ so that for every edge $\{u,v\}\in E(G)$, we have that $w(\{u,v\})$ is the distance from $\varphi(u)$ to $\varphi(v).$ Similarly define a {\it dot product $G$-configuration} where we replace distance by dot product. 

\begin{definition}\label{main}
Given any large, finite set of $n$ points in $\mathbb R^2,$ and a graph $G$, define the maximum number of distance $G$-configurations that can exist in $E$ to be $f(G; n).$ Similarly, define the maximum number of dot product $G$-configurations that can exist in the set to be $g(G; n).$
\end{definition}

With this definition in tow, we now focus on giving estimates on the number of $G$-configurations for various families of edge-weighted graphs. Also note that, in many cases, the particular edge weights will not affect the results, only whether or not an edge is specified in the graph. In light of this, we will not explicitly mention the weights unless they are pertinent. We now name several special families of graphs, and translate known bounds or give new bounds on the configurations they define.

A {\it path} $P_k$ is a graph consisting of $k-1$ distinct edges where each edge shares a vertex with the next for a total of $k$ distinct vertices. Theorems \ref{FKthm}, \ref{2chainzgeneral}, \ref{sharpFKthm}, and \ref{sharpDpChains} can be restated using this notation by expressing upper or lower bounds on $f(P_{k+1};n)$ or $g(P_{k+1}; n).$

\subsubsection{Specific families in two dimensions}

A {\it cycle} $C_k$ is formed by adding an edge to $P_k$ so that the first and last vertices in the path are adjacent. Because any $C_k$ contains a $P_k$, the upper bounds on $f(P_{k+1};n)$ and $g(P_{k+1}; n)$ imply the same upper bounds on $f(C_{k+1};n)$ and $g(C_{k+1}; n),$ respectively. We pause for a moment to focus on the triangle, $C_3.$ In the case of distances, the best bounds that we have for $f(C_3;n)$ follow from the fact any pair of distinct triangles can share no more than one pair of vertices. So for each pair of points $(x, y)$ separated by a given distance $\alpha,$ we have at most 4 choices for $z$ such that $(x,y,z)$ or $(y,x,z)$ form a given triangle (distance $C_3$-configuration). Combining this with the trivial lower bound gives is
\begin{equation}\label{c3}
n\leq f(C_3;n)\lesssim f(P_1; n)\lesssim n^\frac{4}{3}.
\end{equation}
See \cite{BMP}, \cite{DFGMMPS}, and the references contained therein for more details. In the case of dot products, we could potentially fix a pair of points $(x_1,x_3)$ and have as many as $n$ choices for $x_2$ to make a given dot product $C_3$-configuration. However, we can recover the same upper bound with more careful analysis.

\begin{theorem}\label{dpTriangle}
Given any large finite natural number $n,$ we have that
$$n \lesssim g(C_3;n)\lesssim n^\frac{4}{3}.$$
\end{theorem}

A {\it $k$-star} is a graph on $k+1$ vertices where one vertex is adjacent to the other $k$, but none of the other vertices are adjacent to one another.
\begin{theorem}\label{kStar}
Given any large finite natural number $n,$ and an edge-weighted $k$-star $G$, we have that
$$f(G;n)\approx g(G;n)\approx n^k.$$
\end{theorem}

A {\it tree} is a graph with no cycle as a subgraph. In a tree, a vertex with degree one is called a {\it leaf}. A {\it binary tree} is a tree where every vertex has degree one, two, or three. Often times, it is useful to designate one vertex the {\it root} of a tree. Such a tree is then called {\it rooted}. The neighbors of the root are called its {\it children}, and the root is called the {\it parent} of its children. Any neighbors that these vertices have, excluding the root, will be called their children, and so on. The number of edges the path from the root to any vertex is called the {\it height} of that vertex. A {\it perfect binary tree of height $h$} has exactly $2^h$ leaves, each of height $h$. This definition forces every vertex, except for the leaves, to have exactly two children. A $c$-ary tree is a generalization of a binary tree where each parent can have up to $c$ children. So when $c=2,$ this is a binary tree. Let $T_{c,h}$ denote the perfect $c$-ary tree of height $h$.

\begin{theorem}\label{trees}
Given any large finite natural number $n,$ and $c\geq 2,$ we have that
$$f(T_{c,h};n)\approx n^{c^h}.$$
\end{theorem}

We can get a similar result for dot products, but there are technical obstructions to it holding in full generality. So we record here a special case that is relatively easy to state, and explore the more general bound in another estimate that we present later. Let $T_{c,h,\alpha}$ denote an edge-weighted perfect $c$-ary tree of height $h$ whose weights are all the same $\alpha\neq 0.$

\begin{theorem}\label{treesSameDP}
Given any large finite natural number $n,$ and $\alpha\neq 0,$ we have that
$$g(T_{c,h,\alpha};n)\approx n^{c^h}.$$
\end{theorem}

\subsubsection{General families in two dimensions}\label{generalGraphSection}

Where the previous results have been based on rather specific families of graphs, we now move on to some more general families of graphs. We can always get an upper bound on the number of $G$-configurations by decomposing $G$ into subgraphs $H_j$ such that $V(G)=\cup_j V(H_j)$, and taking the product of the upper bounds on the various components. This is in effect ignoring any restrictions imposed on the original configuration by edges between the $H_j.$ That is, we are losing the information provided by the weights of edges in the set $(E(G)\setminus \cup_j E(H_j)).$ In the most extreme case, we would just set each of the $H_j$ to be a single vertex from $G$, and we would get the trivial upper bound of $n^{|V(G)|}.$

To this end, we introduce a definition and illustrate its use to help bridge the gap between the special cases in the previous subsection and more general graphs that could arise. Given any graph $G$, we say that $G$ can be {\it $(s,t)$-covered} if every vertex can be included in a set of $s$ disjoint copies of $P_2$ and $t$ disjoint copies of $P_3.$ Note that if a graph can be $(s, t)$-covered for a pair of parameters $s$ and $t,$ then $|V(G)|=2s+3t,$ as each $P_2$ has two vertices, each $P_3$ has three, and these subgraphs are disjoint.

\begin{minipage}{6in}
\includegraphics[scale=.75]{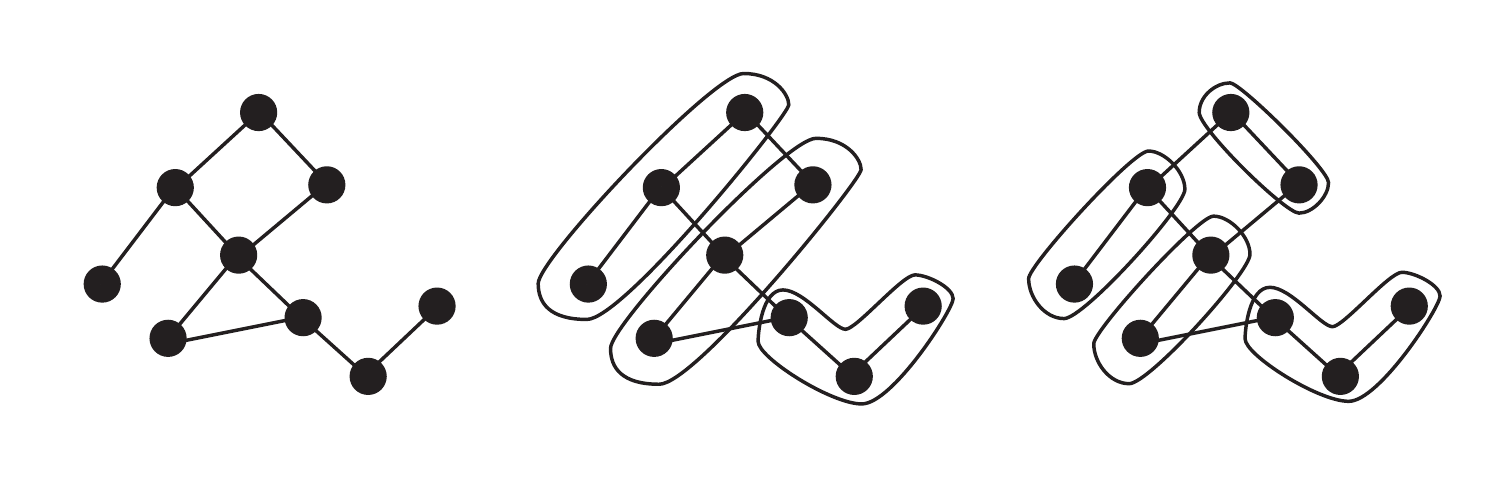}\label{stCover}\\
{\bf Figure 1:} This is a graph that can be $(0,3)$-covered and $(3, 1)$-covered.
\end{minipage}

\begin{theorem}\label{graphCoverThm}
Given any large finite natural number $n,$ and an edge-weighted graph $G$, that can be $(s, t)$-covered, we have that
$$f(G;n)\lesssim n^{\frac{4}{3}s+2t}\text{ and }~g(G;n)\lesssim n^{\frac{4}{3}s+2t}.$$
\end{theorem}

As a simple corollary of Theorem \ref{graphCoverThm}, we get many nontrivial upper bounds for arbitrary $G$-configurations by applying the previous result to a subgraph of $G$ that can be $(s, t)$-covered, and then crudely bounding the number of choices for the remaining vertices by $n$.

\begin{corollary}\label{generalCrude}
Given any large finite natural number $n,$ and an edge-weighted graph $G$, with a subgraph $H$ that can be $(s, t)$-covered, we have that
$$f(G;n)\lesssim n^{\frac{4}{3}s+2t+|V(G)\setminus V(H)|}\text{ and }~g(G;n)\lesssim n^{\frac{4}{3}s+2t+|V(G)\setminus V(H)|}.$$
\end{corollary}

This corollary allows us to prove a dot product version of Theorem \ref{trees}. Here we only state and prove the version for binary trees.

\begin{corollary}\label{dpTrees}
Given any large finite natural number $n,$ 
$$g(T_{2, h};n)\lesssim \begin{cases}
n^{\frac{2}{3}\left(2^{h+1}-1\right)} & \text{ if } h \text{ odd}, \\
n^{1+\frac{4}{3}\left(2^h-1\right)} & \text{ if } h \text{ even}.
\end{cases}$$.
\end{corollary}

\subsubsection{Three dimensional results}

We introduce the notation $f_3(G;n)$, which is the analog of $f(G;n)$ for point sets in $\mathbb R^3.$ We have the following three-dimensional analog of Theorem \ref{kStar}.
\begin{theorem}\label{kStar3}
Given any large finite natural number $n,$ and an edge-weighted $k$-star $G$, we have that
$$f_3(G;n)\approx n^k.$$
\end{theorem}

The situation is a bit different for the three-dimensional version of Theorem \ref{trees}. When $c\geq 3,$ we get similar result, again by restricting that all of the edge weights are the same $\alpha\neq 0.$
 
\begin{theorem}\label{3trees}
Given any large finite natural number $n,$ and a $c\geq 3,$ we have that
$$f_3(T_{c, h,\alpha};n)\approx n^{c^h}.$$
\end{theorem}

However, in the special case that $c=2$, we can no longer be sure the the upper bound matches the lower bound, so we have the following more precise yet technical estimate.

\begin{theorem}\label{techTree}
Given a large natural number $n$, an $\alpha\neq 0,$ and any $\epsilon>0,$
$$f(T_{c, h,\alpha};n)\lesssim n^{c^h\left(1+\frac{1}{c^2-1}\right)+\frac{(-1)^h}{2(c+1)}+\frac{1}{2(1-c)}+\epsilon}$$
\end{theorem}

Applying this in the case that $c=2$ implies the following.

\begin{corollary}\label{binary3Trees}
Given any large finite natural number $n,$ we have that
$$n^{2^h}\lesssim f_3(T_{2, h,\alpha};n)\lesssim n^{\frac{1}{6}\left(2^{h+3}-3+(-1)^h\right)} .$$
\end{corollary}

Finally, in the absence of any other discernible graph structures that we can identify, we have a three-dimensional companion to Theorem \ref{generalCrude}. A {\it matching} in a graph is a set of edges that share no vertices. For a graph $G$, let $m_G$ denote the maximum number of edges in a matching of $G$, and let $r_G$ be the number of remaining vertices. Specifically,
$$m_G:=\max\{E(H):H\text{ is a matching in } G\} \text{ and } r_G:=|V(G)-2m_G|.$$ 

\begin{theorem}\label{lego}
Given any large finite natural number $n,$ any $\epsilon >0,$ and an edge-weighted graph $G$, we have 
$$f_3(G;n)\lesssim n^{\frac{295}{197}m_G+r_G+\epsilon}.$$
\end{theorem}

We conclude the statement of our main results with a concrete comparison. Let us fix our attention on a perfect binary tree of height three whose edge weights are all the same $\alpha\neq 0.$ Theorem \ref{binary3Trees} gives an upper bound of $f_3(T_{2,3,\alpha})\lesssim n^{\frac{14}{3}+\epsilon},$ for any $\epsilon>0.$ If we apply Theorem \ref{lego}, we first look for the largest matching we can find in $T_{2,3,\alpha}.$ The best we can do is $m_{T_{2,3,\alpha}}=2$, which gives us three vertices left over, or $r_{T_{2,3,\alpha}}=3.$ So we get an upper bound of $f_3(T_{2,3,\alpha})\lesssim n^{\frac{590}{197}+3+\epsilon},$ for any $\epsilon >0,$ which is considerably worse than the estimate given by Theorem \ref{binary3Trees}.

\subsection{Organization of this paper}
In Section 2, we flesh out some of the tools that we will use to prove the main results. To this end, we outline proofs of some known results mentioned above. Section 3 has the proofs of the new results.

\section{Preliminaries}\label{prelim}

Distance configurations can be explored by looking at the incidences of points and families of circles centered at those points. The analogous geometric objects for dot product configurations are special families of lines.

\subsection{The $\alpha$-line for a point $p$}

Given a point $p\in \mathbb R^2\setminus\{(0,0)\}$, the set of points in the plane that have dot product $\alpha$ with $p$ is a line we call the {\it $\alpha$-line} of $p$, $\ell_\alpha(p).$ We call any line through the origin a {\it radial} line. A quick calculation using the definition of the dot product gives us that for any $p\in \mathbb R^2\setminus\{(0,0)\}$ and any real number $\alpha,$ we have that $\ell_\alpha(p)$ will be perpendicular to the unique radial line through $p$. We include the following lemmata from \cite{BS}.

\begin{lemma}\label{singlePoint}
If $p$ and $r$ are two points in $\mathbb{R}^2\setminus\{(0,0)\}$ that do not lie on the same radial line, and $\alpha, \beta \in \mathbb{R},$ then there exists exactly one point $q \in\mathbb{R}^2\setminus\{(0,0\}$ such that $p \cdot q = \alpha$ and $q \cdot r = \beta.$
\end{lemma}

Given a point $p\in \mathbb R^2,$ let $|p|$ denote the distance from $p$ to the origin.

\begin{lemma}\label{setB}
If $p$ and $r$ are two points in $\mathbb{R}^2\setminus\{(0,0)\}$ such that $\ell_\alpha(p)$ coincides with $\ell_\beta(r),$ then $p$ and $r$ must lie on the same radial line through the origin, and $|p|/|r|=\beta/\alpha.$
\end{lemma}

\begin{minipage}{5in}
\includegraphics[scale=.75]{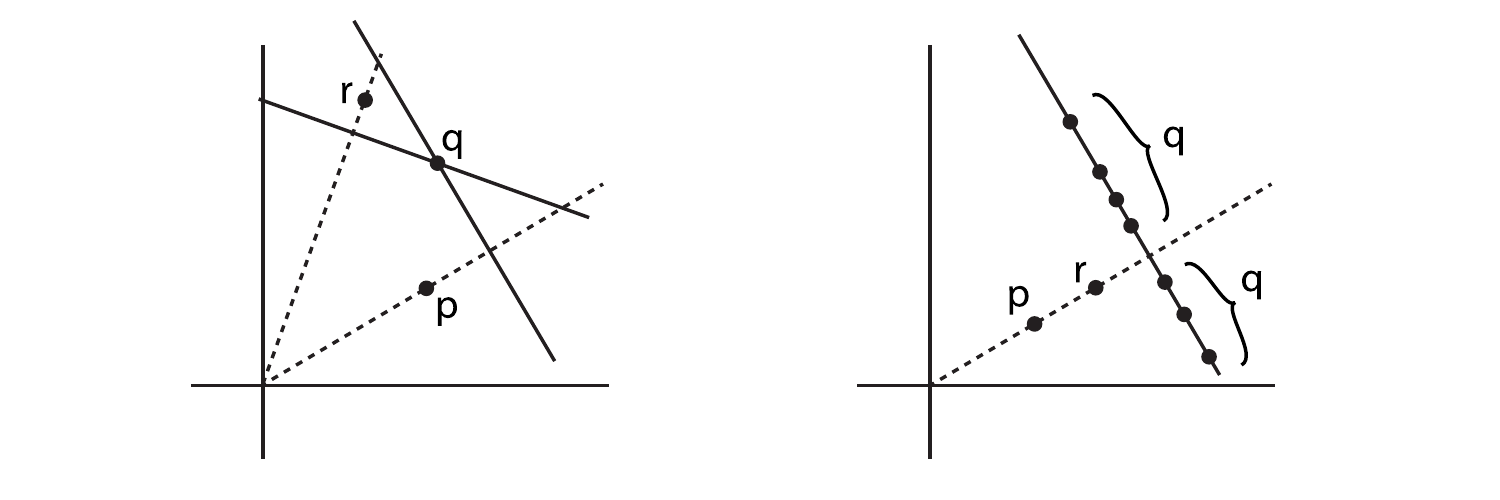}\label{alphaLines}\\
{\bf Figure 2:} Here, we are looking for triples of points $(p,q,r)$ such that $p\cdot q=\alpha$ and $q\cdot r=\beta.$ On the left, we have points $p$ and $r$ on distinct radial lines, so $q$ is the unique intersection of $\ell_\alpha(p)$ and $\ell_\beta(r)$ for positive numbers $\alpha$ and $\beta.$ On the right, we see that $p$ and $r$ are on the same radial line, so for some specific choices of positive real numbers $\alpha$ and $\beta$, we can have $\ell_\alpha(p)$ coincide with $\ell_\beta(r),$ giving potentially many choices for points $q$.\\\\
\end{minipage}

\subsection{Sketch of proof of Theorem \ref{hinges}}
We take a moment to give an outline of the proof of Theorem \ref{hinges} to illustrate the basic ideas used in the main results below. This highlights some of the fundamental similarities and differences between distances and dot products.

In the case of distances, Theorem \ref{hinges} follows from the observation that there are about $n^2$ pairs of points $(x_1,x_3)$ and for each such pair, there at most two intersections of the circles of radii $\alpha_1$ and $\alpha_2$ centered at $x_1$ and $x_3$, respectively. These two intersections are the only possible locations for $x_2$ if $(x_1, x_2, x_3)$ is to be a distance hinge of type $(\alpha_1, \alpha_2).$

The companion result for dot product hinges is similar, but a bit more involved. We follow the same program as with distances, but for each pair of points, $(x_1,x_3),$ we apply either Lemma \ref{singlePoint} or Lemma \ref{setB}, depending on whether or not $\ell_{\alpha_1}(x_1)$ coincides with $\ell_{\alpha_2}(x_3).$ In the first case, everything follows as with distances, but in the second case, we have to count a bit more carefully. See \cite{BS} for more details.

\subsection{The Szemer\'edi-Trotter Theorem}

The celebrated Szemer\'edi-Trotter Theorem from \cite{ST83} is a key component of many of the dot product proofs. We also note that this result provided the foundation for Theorem \ref{SST}, which is heavily leaned on in the study of distance configurations.

\begin{theorem}\label{ST}
Given $n$ points and $m$ lines in the plane, the number of point-line pairs, such that the point lies on the line is
$$\lesssim\left(n^{\frac{2}{3}}m^{\frac{2}{3}}+n+m\right).$$
\end{theorem}

Theorem \ref{STdp} follows from this by considering the set of $n$ points in the plane, and the set of their $\alpha$-lines. Lemma \ref{setB} then guarantees that there will be $n$ unique such lines, and the result follows.

\section{Proofs}\label{proofs}

\subsection{Proof of Theorem \ref{dpTriangle}}
The lower bound follows using a geometric progression along any line through the origin. We give an explicit construction for points along the $x$-axis.
Let the set of points be $\{(ag^j,0): j=1\dots n\}$, for $a, g >0.$ We then have $\gtrsim n$ triples of points of the form $(x_1, x_2, x_3)=((ag^{j-1},0), (ag^j,0), (ag^{j+1},0)),$
with the property that $x_1\cdot x_2 = a^2g^{2j-1},$ $x_2\cdot x_3 = a^2g^{2j+1},$ and $x_3\cdot x_1 = a^2g^{2j},$ and the lower bound is achieved.

For the upper bound, we want to get a bound on the number of triples of points of the form $(x_1, x_2, x_3)$ such that $x_1\cdot x_2 = \alpha_1,$ $x_2\cdot x_3 = \alpha_2,$ and $x_3\cdot x_1 = \alpha_3.$ We start by considering every pair of points from the set that could possibly be $x_1$ and $x_3$ from such a triple. Theorem \ref{STdp}, we know that there can be no more than $\lesssim n^\frac{4}{3}$ point pairs that have dot product $\alpha_3,$ so we know that there can be no more than $\lesssim n^\frac{4}{3}$ choices for the pair $(x_1,x_3).$

For each such pair of points, either they do not lie on the same radial line or they do lie on the same radial line. If they do not lie on the same radial line, then Lemma \ref{singlePoint} guarantees that $\ell_{\alpha_1}(x_1)$ intersects $\ell_{\alpha_2}(x_3)$ in exactly one point, and this point is the only potential candidate for an $x_2$ that would make $(x_1, x_2, x_3)$ a dot product $C_3$-configuration. Therefore, there are at most $n^\frac{4}{3}$ dot product $C_3$-configurations with $x_1$ and $x_3$ on distinct radial lines.

We now turn our attention to the pairs of points $(x_1, x_3)$ that satisfy $x_3\cdot x_1 = \alpha_3$ and have both points on the same radial line. Because $x_1$ and $x_3$ share the same radial line, we know that $\ell_{\alpha_1}(x_1)$ and $\ell_{\alpha_2}(x_3)$ will have the same slope, as both of these lines must be perpendicular to the radial line containing $x_1$ and $x_3$. If $\ell_{\alpha_1}(x_1)$ does not coincide with $\ell_{\alpha_2}(x_3),$ then they must not intersect, as they have the same slope. Therefore, there is no possible location for a point $x_2$ that would make $(x_1, x_2, x_3)$ a dot product $C_3$-configuration.
\subsubsection{Coincident $\alpha$-lines}\label{coincident}
So the last (and most involved) possibility to check is when $x_1$ and  $x_3$ satisfy $x_3\cdot x_1 = \alpha_3$ and we have that $\ell_{\alpha_1}(x_1)$ coincides with $\ell_{\alpha_2}(x_3)$. In this case, we appeal to Lemma \ref{setB}, which tells us that
$$
\frac{|x_1|}{|x_3|} = \frac{\alpha_1}{\alpha_2}.
$$
\begin{equation}\label{x1}
|x_1|=|x_3|\frac{\alpha_1}{\alpha_2}.
\end{equation}
But recall that $x_3\cdot x_1 = |x_3||x_1|\cos \theta,$ where $\theta$ is the angle made by $x_1$, the origin, and $x_3.$ But $x_1$ and $x_3$ are on the same radial line, so $\cos \theta = \pm 1.$ Combining this with the fact that $x_3\cdot x_1 = \alpha_3$, gives us
$$\pm\alpha _3=|x_3||x_1|.$$
Combining this with \eqref{x1} yields
$$\pm\alpha_3=|x_3|^2\frac{\alpha_1}{\alpha_2},$$
which means that $x_3$ must live on a circle of radius $\sqrt{\alpha_3\alpha_2/\alpha_1}$ centered at the origin. Call this circle $\mathcal C.$ Recall that since $\ell_{\alpha_1}(x_1)$ and $\ell_{\alpha_2}(x_3)$ coincide, any point on this line will be a suitable choice for $x_2$ to make $(x_1, x_2, x_3)$ a dot product $C_3$-configuration. We now pause to state a small fact about antipodal points on a circle, whose technical proof we delay.

\begin{claim}\label{antipodal}
If $p$ and $p'$ are antipodal points on a circle centered at the origin, then for any nonzero $\alpha,$ we have that $\ell_{\alpha}(p)$ is distinct from $\ell_{\alpha}(p')$.  
\end{claim}

So we will finish this case by estimating the number of incidences of points from our set (putative choices for $x_2$) and lines of form $\ell_{\alpha_2}(x_3)$ for choices of $x_3$ on the circle $\mathcal C.$ Notice that each distinct choice of $x_3$ must give rise to a different line of the form $\ell_{\alpha_2}(x_3)$, because each point $x_3\in \mathcal C$ will have lie on a different radial line, except for its antipodal point $x_3',$ but in that case, we can apply Claim \ref{antipodal} to guarantee that $\ell_{\alpha_2}(x_3)$ will not coincide with $\ell_{\alpha_2}(x_3')$ for nonzero choices of $\alpha_2.$ Therefore, our final count amounts to estimating the incidences between $\lesssim n$ points and $\lesssim n$ distinct lines in the plane. By appealing to Theorem \ref{ST}, we see that there can be no more than $n^\frac{4}{3}$ such incidences, completing the proof.

\subsubsection{Proof of Claim \ref{antipodal}}
\begin{proof}
By way of contradiction, suppose $\ell_{\alpha}(p)$ coincides with $\ell_{\alpha}(p').$ Because $p$ and $p'$ are antipodal points on a circle centered at the origin, they must lie on the same radial line. Let $q$ be the intersection of $\ell_{\alpha}(p)$ and the radial line through $p$ and $p'.$ Let $\theta$ be the angle determined by $p,$ the origin, and $q$. Define $\phi$ to be the angle between $p',$ the origin, and $q$. Because $p$, $p'$, $q$, and the origin all lie on a line, with $p$ and $p'$ on opposite sides of the origin, we must have that $\{ \theta, \phi \} = \{0, \pi \}.$ However, by our assumptions and the definition of dot product, we have that
$$|p||q|\cos \theta = p\cdot q = \alpha = p'\cdot q = |p'||q|\cos \phi.$$
Because $p$ and $p'$ lie on the same circle, we know that $|p|=|p'|.$ So comparing the left and right hand sides of the above equation yields $\cos \theta = \cos \phi,$ which implies that $1=-1,$ which is a contradiction.
\end{proof}

\subsection{Proof of Theorem \ref{kStar}}
Start by labeling the vertices of $G$ as $v_1, \dots, v_{k+1},$ where each of the vertices $v_1$ through $v_k$ has degree 1, and the vertex $v_{k+1}$ has degree $k$.

Now we show that $f(G;n)\lesssim n^k.$ To see this, we apply Theorem \ref{FKthm} for 2-chains to see that there are no more than $n^2$ triples of points of the form $(x_1,x_{k+1}, x_2)$ such that $|x_1-x_{k+1}|=w(v_1,v_{k+1})$ and $|x_2-x_{k+1}|=w(v_2,v_{k+1})$. For each such triple, we then have no more than $n$ choices for each of the other $x_j,$ for a total of $n^k.$

We then prove that $f(G;n)\gtrsim n^k.$ Start by assuming the origin is in our set. Then arrange about $n/k$ points on $k$ circles centered at the origin, whose radii are the edge weights, $w(v_j,v_{k+1}).$ If any weights repeat, just put more points on the circle of the appropriate radius until we have placed about $n$ points.

The same results follow for the upper and lower bounds of $g(G;n),$ except by appealing to Theorem \ref{hinges} for the upper bound, and with appropriate $\alpha$-lines in place of circles for the lower bound.

\subsection{Proofs of Theorem \ref{trees} and Theorem \ref{treesSameDP}}

We start by proving Theorem \ref{trees}, and later describe how to modify the proof to prove Theorem \ref{treesSameDP}. Call the root of $G$ $v_1.$ Let the children of $v_1$ be called $v_2,v_3, \dots v_{c^2-1}.$ Continuing in this way, we get that the vertices at height $k$ will have indices $i$ satisfying $c^k\leq i <c^{k+1}.$ Because $G$ is perfect and of height $h$, we will have a total of $c^{h+1}-1$ vertices. Now, we have $c^h$ leaves, and $n$ choices for each of them, giving us a total of $n^{2^h}$ choices for set of leaves. Every pair of leaves sharing a parent in the graph determine exactly two possibilities for points from the set that could serve as parents, just as in the proof of Theorem \ref{hinges}. Since there will be $c^{h-1}$ such relationships, the total number of distance $G$-configurations is no more than $c^{h-1}n^{c^h},$ as claimed.

To see the lower bound, we will construct a point set, $E$, with $\gtrsim n^{c^k}$ instances of a given $G$-configuration. Consider $H$, the induced subgraph of $G$ formed by removing the leaves from $G$. Explicitly construct a single instance of an $H$-configuration by picking a point to correspond to the root, then picking points at the prescribed distances from the root corresponding to the first generation of children, and so on, until we have a single instance of the prescribed $H$-configuration. We will finish by processing each leaf of $H$ in the following way. Suppose $v$ is a leaf in $H$. Then in $G$, $v$ has $c$ children. Let $w_j$ be the weight of the edge between $v$ and its $j$th child. Now, for each child, arrange $c^{-h}n$ points on the sphere of radius $w_j$ centered at the point of $E$ corresponding to the vertex $v$. Do this for every leaf of $H$. Notice that by construction, we have $c^{-h}n\approx n$ choices of points to correspond to every leaf of $G$, and moreover, these choices are independent.

To prove Theorem \ref{treesSameDP}, we follow the same logic, except with $\alpha$-lines in place of circles. The restriction that all of the edge weights in $G$ are the same nonzero value $\alpha$ prevents the existence of potentially pathological sub-configurations (such as in the final case of the proof of Theorem \ref{dpTriangle}) where the points corresponding to two leaf vertices might not uniquely determine the point corresponding to their parent vertex.

\subsection{Proof of Theorem \ref{graphCoverThm}}

Because $G$ is $(s, t)$-covered, we know that for each occurrence of the given distance $G$-configuration, we must have at least $s$ pairs of points separated by fixed distances determined by the weights on their corresponding edges, and at least $t$ triples of points that form hinges with distance relations also defined by the appropriate edge weights. Theorem \ref{SST} guarantees that there can be no more than $n^\frac{4}{3}$ pairs of points separated by any given distance, and Theorem \ref{hinges} tells us that there can be no more than $n^2$ occurrences of any hinge. Putting these together, we get
$$f(G; n)\lesssim \left(n^\frac{4}{3}\right)^s \left(n^2\right)^t,$$
as claimed. To prove the analogous result for dot products, we follow the same logic, but appeal to Theorem \ref{STdp} in place of Theorem \ref{SST} above.

\subsection{Proof of Corollary \ref{dpTrees}}
We begin with the case that $h$ is odd. If $h=1,$ our graph is just $P_2,$ and this amounts to the hinge bound, Theorem \ref{hinges}. If $h=3$, we cover the root and its children by a $P_2,$ and the other two generations are covered by four copies of $P_2.$ In general, we see that the vertices corresponding to every other generation of the tree can be decomposed into copies of $P_2,$ by greedily separating out hinges from every other generation. Because $h$ is odd, we know that there will be an even number of generations, so no points will be left out. As there are $2^{h+1}-1$ vertices in $G$, and they are all broken up into triples, the total number of triples is $\frac{1}{3}\left(2^{h+1}-1\right)$. Theorem \ref{hinges} has guaranteed that no such triple can occur more than $n^2$ times, for a total upper bound of
$$g(G; n)\lesssim n^{\frac{2}{3}\left(2^{h+1}-1\right)}.$$
In the case that $h$ is even, note that $G$ can be decomposed into a root and two perfect binary trees of height $h-1$. So we apply the odd height estimate to each of these trees, and have no more than $n$ choices for the root, yielding
$$g(G; n)\lesssim n^{\frac{2}{3}\left(2^h-1\right)}n^{\frac{2}{3}\left(2^h-1\right)}n.$$

\begin{minipage}{6in}
\includegraphics[scale=.75]{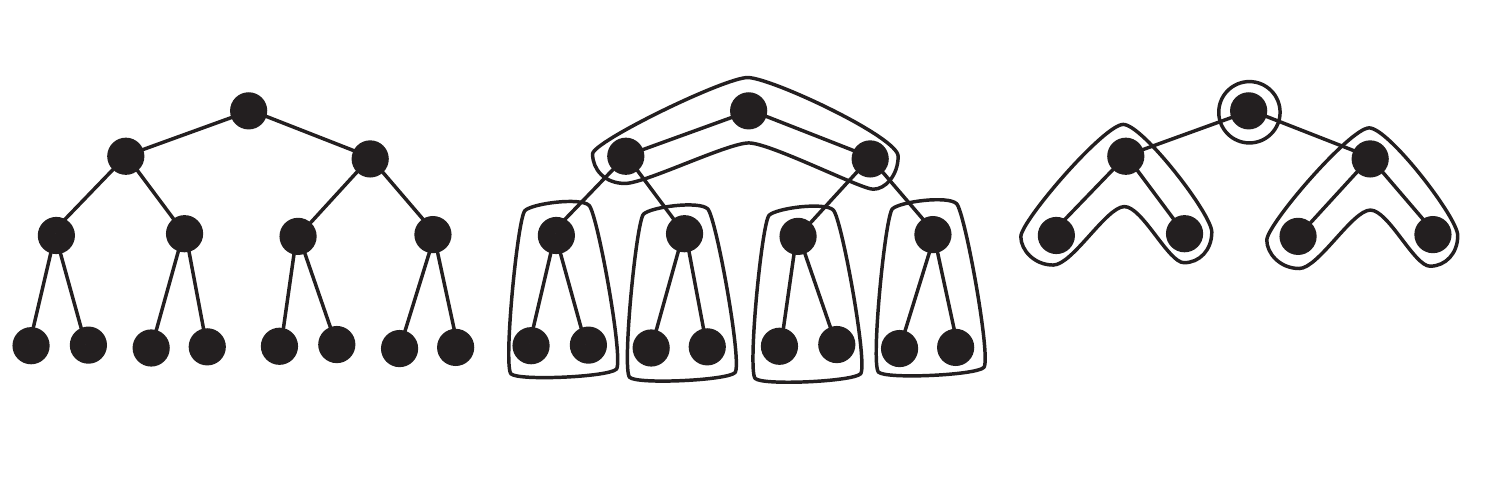}\label{binaryTree}\\
{\bf Figure 3:} This is a perfect binary tree of height 3, the same tree with a\\ $(0,5)$-covering shown, and a perfect binary tree of height 2 with a\\ $(0,2)$-covering on all but its root, which must be considered separately.
\end{minipage}

\subsection{Proof of Theorem \ref{3trees}}

The lower bounds follow from a similar construction to the one in the proof of Theorem \ref{trees}, except with spheres in place of circles. The proof of the upper bounds also follow by similar reasoning, except that instead of the intersections of two circles determining constantly many points, we use that three spheres of the same radius will intersect in constantly many points.

\subsection{Proof of Theorem \ref{techTree}}
The basic idea is to decompose the tree into subgraphs that can be estimated separately. We begin by estimating the longest chain, which is a path from one leaf, through the root, and to another leaf. Therefore, this path will be $P_{2h+1},$ which is a distance $2h$-chain. We will call this path the {\it outer chain}. Next, there will be $c-2$ perfect $c$-ary trees of height $h-1$ whose roots were children of the original root. This is because the root of the big tree has $c$ children, but two of them were accounted for in the outer chain.

Notice that each of these two children of the root will have one child in the outer chain, and their other $c-1$ children (grand-children of the root) as yet unaccounted for. This gives us a total of $2c-2$ children that are not yet counted in this generation. Each of these, being at height two from the root, will then be the root of their own perfect $c$-ary trees of height $h-2.$ This pattern continues with the next children of vertices in the outer chain that are as yet unaccounted for. There will again be $c-1$ children from each of the outer chain vertices, and these children will again be the roots of their own perfect $c$-ary trees of successively lower heights. Putting all of these pieces together, we get the following recursive relationship.
\begin{equation}\label{bigTree}
f(T_{c, h,\alpha};n)\lesssim f(P_{2h+1}; n) f(T_{c, h-1,\alpha}; n)^{c-2} \prod_{j=0}^{h-2}f(T_{c, j,\alpha})^{2c-2}.
\end{equation}

To simplify the analysis, we will let $a_{c, h}$ be the exponent of $n$ from the best estimate we have for $f(T_{c, h,\alpha};n)$. That is, $$a_{c, h} = \min\{x: \text{ we can show } f(T_{c, h,\alpha};n)\lesssim n^x\}.$$ We will work with $a_{c, h},$ acknowledging that it might not correspond to the best possible estimate for $f(T_{c, h,\alpha};n),$ but it will give us an upper bound. Continuing, we translate what we know about each of the components of \eqref{bigTree} into expressions related to $a_{c, h}.$

Recall that a $P_{2h+1}$-configuration is just a $2h$-chain, so by Theorem \ref{FKthm3}, we know that $f(P_{2h+1};n)\lesssim n^{\frac{2h}{2}+1+\epsilon},$ for any $\epsilon >0.$ So this term will give us a factor of $n^{h+\epsilon}.$ We also notice that $T_{c, 0,\alpha}$ is just a single point, so $f(T_{c, 0,\alpha};n)=n.$ This gives us that $a_{c,0}=1.$ Similarly, we can see that $T_{c,1,\alpha}$ is a $c$-star. So Theorem \ref{kStar3} tells us that $f(T_{c, 1,\alpha};n)\approx n^c.$ This gives us that $a_{c,1}=c.$ Now we can rewrite \eqref{bigTree} in terms of $a_{c, h}.$ We get the following recurrence relation, with two initial conditions.
\begin{equation}\label{aTree}
a_{c, h}=(h+\epsilon)+(c-2)a_{c, h-1}+\sum_{j=0}^{h-2}(2c-2)a_{c, j}; ~~a_{c,0}=1; ~~a_{c,1}=c.
\end{equation}
This recurrence relation has the solution
$$a_{c,h}=c^h\left(1+\frac{1}{c^2-1}\right)+\frac{(-1)^h}{2(c+1)}+\frac{1}{2(1-c)}+\epsilon.$$
This implies the desired result.

\begin{minipage}{6in}
\includegraphics[scale=.75]{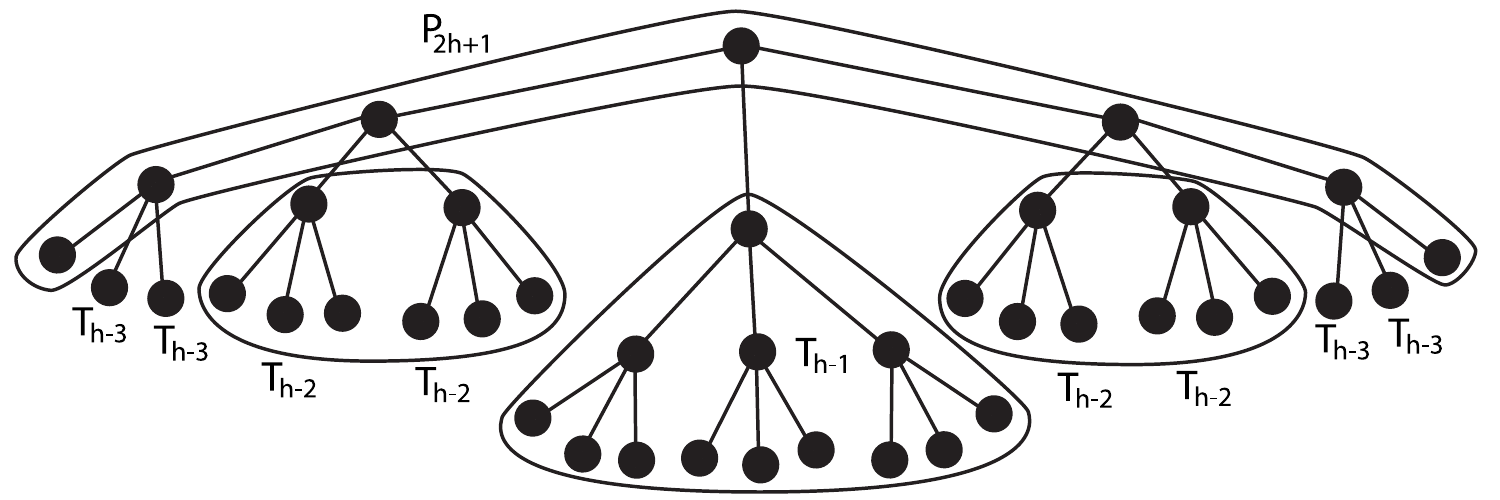}\label{binaryTree}\\
{\bf Figure 4:} This is the root and next three generations of a $T_{3,3}$, with\\
the various components labeled as in the proof of Theorem \ref{binary3Trees}.
\end{minipage}

\subsection{Proof of Theorem \ref{lego}}

This proof is similar to that of Theorem \ref{graphCoverThm}. We find a maximal matching of $G$ consisting of $m_G$ edges. Each of these edges can have no more than $n^{\frac{295}{197}+\epsilon}$ pairs of points that can represent them. There are no more than $n$ choices for any of the $r_G$ vertices of $G$ that have yet to be accounted for. Putting this all together yields the desired result.

\end{document}